\documentclass[preprint,review,11pt]{elsarticle}

\usepackage{lineno,hyperref}
\modulolinenumbers[5]

\journal{Journal of \LaTeX\ Templates}

\everydisplay{\mathsurround 0pt} 

\usepackage{amsmath,amssymb,amsthm,mathtools}

\newcommand{\R}{\mathbb{R}}

\newtheorem{theorem}{Theorem}
\newtheorem{corollary}{Corollary}
\usepackage{setspace}
\begin{document}

\begin{frontmatter}
\title{A general measure of the impact of priors in Bayesian statistics via Stein's Method}

\author{Fatemeh Ghaderinezhad}
\ead{fatemeh.ghaderinezhad@ugent.be}

\author{Christophe Ley\corref{cor2}}
\ead{christophe.ley@ugent.be}
\cortext[cor2]{Corresponding author}
\address{Ghent University, Department of Applied Mathematics, Computer Science and Statistics\\Krijgslaan 281, S9, Campus Sterre, 9000 Ghent, Belgium}


\begin{abstract}
\indent We propose a measure of the impact of any two choices of prior distributions by quantifying the Wasserstein distance between the respective resulting posterior distributions at any fixed sample size. We illustrate this measure on the normal, Binomial and Poisson models.
%
%

\end{abstract}


\end{frontmatter}

\section{Introduction}

A key question in Bayesian analysis is the effect of the prior on the posterior, and how we can measure this effect. Will the posterior distributions derived with distinct priors become very similar if more and more data are gathered? It has been proved formally in~\citet{DF1} and \citet{DF2} that, under certain regularity conditions, the impact of the prior is waning as the sample size increases. From a practical viewpoint it is more important to know what happens at finite sample size $n$. Recently, \citet{LRS17a} have provided a partial answer to this question by investigating the Wasserstein distance between the posterior distribution based on a given prior of interest and the  no-prior data-only based posterior, which allows detecting at fixed sample size $n$ the effect of the prior of interest. This distance being mostly impossible to calculate exactly, they have provided sharp upper and lower bounds. This work is, to the best of the authors' knowledge, the first to quantify at any sample size the prior effect. However, it strongly relies on the assumption that one prior is the flat uniform prior (or data-only prior), and hence it does not allow a direct comparison between two general priors.

Our aim in the present paper is to fill this gap and to extend the methodology of \citet{LRS17a} to the setting where we can quantify the distributional distance between posteriors resulting from any two priors. As in \citet{LRS17a}, we opt here  for the Wasserstein-1 distance defined as
$$
d_{\mathcal W}(P_1,P_2)=\sup_{h\in\mathcal{H}}|{\rm E}[h(X_1)]-{\rm E}[h(X_2)]|
$$
for $X_1$ and $X_2$ random variables with respective distribution functions $P_1$ and $P_2$, and where $\mathcal{H}$ stands for the class of Lipschitz-1 functions. We shall provide sharp lower and upper bounds on this distance for any two posteriors. For the upper bound case, our approach relies on a variant of  Stein's density approach \citet{Stein2004}, see also \citet{LRS17b} for a recent account. 

The present paper is organized as follows. In Section~\ref{sec2} we provide the notations and terminology used throughout the paper, state and prove our main theorem. Then in Section~\ref{sec3} we illustrate the strength of our new general measure of the difference between two priors by considering several examples of well-known distributions. In each case, we compare the effects of distinct priors on the posterior distribution of the parameter we are interested in. The choices of our priors are motivated by research papers that have discussed various choices of viable priors for a certain parameter.

\section{A general method to compare the effects of two priors}\label{sec2}

We start by fixing our notations. We consider independent and identically distributed observations $X_1,\ldots,X_n$  from a parametric model with parameter of interest  $\theta \in \Theta\subseteq\R$. Our framework allows the common distribution of the $X_i$'s to be continuous as well as discrete. For the sake of simplicity, we shall denote the sampling density or likelihood of $X_1,\ldots,X_n$ by $\ell(x;\theta)$ where $x=(x_1,\ldots,x_n)$ are the observed values. Now consider two distinct (possibly improper) prior densities  $p_1(\theta)$ and  $p_2(\theta)$ for $\theta$.  The two resulting posterior densities for $\theta$ then take on the guise
$$
p_i(\theta;x)=\kappa_i(x)p_i(\theta) \ell(x;\theta), \quad i=1,2,
$$
where  $\kappa_1(x),\kappa_2(x)$ are normalizing constants. We write $(\Theta_1,P_1)$ and $(\Theta_2,P_2)$ the couples of random variables and cumulative distribution functions corresponding respectively to the densities $p_1(\theta;x)$ and $p_2(\theta;x)$. 

Suppose that both posterior distributions have finite means $\mu_1$ and $\mu_2$. For our general result we borrow a quantity from the Stein Method literature, namely the \emph{Stein kernel} $\tau_i$ of $P_i$ defined as
\begin{equation}\label{kernel}
\tau_i(\theta;x)=\frac{1}{p_i(\theta;x)}\int_{a_i}^{\theta}(\mu_i-y)p_i(y;x)dy, \quad i=1,2,
\end{equation}
where $a_i$ is the lower bound of the support $I_i=(a_i,b_i)$ of $p_i$. One can readily see that this function is always positive and vanishes at the boundaries of the support. For more information about Stein kernels, we refer the reader to \citet{LRS17a}; for the present paper no further knowledge than the definition~\eqref{kernel} is necessary.

The key element in our developments is the fact that the densities $p_1(\theta;x)$ and $p_2(\theta;x)$ are \emph{nested}, meaning that one support is included in the other. We here suppose that ${I}_2\subseteq{I}_1$ which allows us to write $p_2(\theta;x)=\frac{\kappa_2(x)}{\kappa_1(x)}\rho(\theta) p_1(\theta;x)$ with 
$$
\rho(\theta)=\frac{p_2(\theta)}{p_1(\theta)}.
$$
With these notations in hand, we obtain the following theorem which is a general way of comparing the two priors $p_1(\theta)$ and $p_2(\theta)$ and follows from Stein's Method for nested densities \citet{LRS17a}.
\begin{theorem}\label{maintheo}
Consider $\mathcal{H}$ the set of Lipschitz-1 functions on $\R$. Assume that $\theta\mapsto\rho(\theta)$ is differentiable on ${I}_2$ and satisfies (i) ${\rm E}[|\Theta_1-\mu_1|\rho(\Theta_1)]<\infty$, (ii) $\left(\rho(\theta)\int_{a_1}^{\theta}(h(y)-{\rm E}[h(\Theta_1)])p_1(y;x)dy\right)'$ is integrable for all $h\in\mathcal{H}$ and (iii) $\lim_{\theta\rightarrow a_2,b_2}\rho(\theta)\int_{a_1}^{\theta}(h(y)-{\rm E}[h(\Theta_1)])p_1(y;x)dy=0$ for all $h\in\mathcal{H}$. Then
\begin{equation}\label{bounds}
|\mu_1-\mu_2|=\frac{|{\rm E}[\tau_1(\Theta_1;x)\rho'(\Theta_1)]|}{{\rm E}[\rho(\Theta_1)]}\leq d_{\mathcal{W}}(P_1,P_2)\leq \frac{{\rm E}[\tau_1(\Theta_1;x)|\rho'(\Theta_1)|]}{{\rm E}[\rho(\Theta_1)]}
\end{equation}
and, if the variance of $\Theta_1$ exists,

\begin{equation}\label{bounds2}
|\mu_1-\mu_2|\leq d_{\mathcal{W}}(P_1,P_2)\leq ||\rho'||_{\infty}\frac{{\rm Var}[\Theta_1]}{{\rm E}[\rho(\Theta_1)]}
\end{equation}
where $||\cdot||_{\infty}$ stands for the infinity norm.
\end{theorem}
\emph{Proof.} Our proof relies on Theorem 3.1 of \citet{LRS17a} about Stein's Method for nested densities. Let us start by establishing the upper bounds in~\eqref{bounds} and~\eqref{bounds2}. The function $\pi_0$ from \citet{LRS17a} coincides in our case with $\frac{\kappa_2(x)}{\kappa_1(x)}\rho(\theta)$, yielding our conditions on $\rho(\theta)$ and the upper bound
$$
{\rm E}[\tau_1(\Theta_1;x)|\rho'(\Theta_1)|]\frac{\kappa_2(x)}{\kappa_1(x)}.
$$
The final form of the upper bound in~\eqref{bounds} is obtained by noticing that 
$$
1=\int_{\Theta}p_2(\theta;x)d\theta=\int_{\Theta}\frac{\kappa_2(x)}{\kappa_1(x)}\rho(\theta) p_1(\theta;x)d\theta=\frac{\kappa_2(x)}{\kappa_1(x)}{\rm E}[\rho(\Theta_1].
$$
The upper bound in~\eqref{bounds2} follows from ${\rm E}[\tau_1(\Theta_1;x)|\rho'(\Theta_1)|]\leq ||\rho'||_{\infty}{\rm E}[\tau_1(\Theta_1;x)]$ and from the property of a Stein kernel to satisfy ${\rm E}[\tau_1(\Theta_1;x)\varphi'(\Theta_1)]={\rm E}[(\Theta_1-\mu_1)\varphi(\Theta_1)]$ for all differentiable function $\varphi$ for which ${\rm E}[(\Theta_1-\mu_1)\varphi(\Theta_1)]$ exists, implying in particular that ${\rm E}[\tau_1(\Theta_1;x)]={\rm Var}[\Theta_1]$. 

\

The lower bound $\frac{|{\rm E}[\tau_1(\Theta_1;x)\rho'(\Theta_1)]|}{{\rm E}[\rho(\Theta_1)]}$  is  obtained along the same lines as the upper bound. It remains to show that it  equals the simplified expression $|\mu_1-\mu_2|$. Using again the above-mentioned Stein kernel property, we get
\begin{eqnarray*}
\frac{|{\rm E}[\tau_1(\Theta_1;x)\rho'(\Theta_1)]|}{{\rm E}[\rho(\Theta_1)]}&=&|{\rm E}[\tau_1(\Theta_1;x)\rho'(\Theta_1)]|\frac{\kappa_2(x)}{\kappa_1(x)}\\
&=&|{\rm E}[(\Theta_1-\mu_1)\rho(\Theta_1)]|\frac{\kappa_2(x)}{\kappa_1(x)}\\
&=&\left|{\rm E}[\Theta_2]-\mu_1{\rm E}[\rho(\Theta_1)]\frac{\kappa_2(x)}{\kappa_1(x)}\right|\\
&=&|\mu_2-\mu_1|.
\end{eqnarray*}
\hfill $\square$

\

We attract the reader's attention to the fact that the lower bound $|\mu_1-\mu_2|$ can be readily obtained by noticing that the identity function $h(x)=x$ belongs to the class of Lipschitz-1 functions. For the comparison of lower and upper bounds it is however relevant to provide as well the equivalent expression of the lower bound in terms of the Stein kernel and $\rho$. Indeed, a close examination of both upper and lower bounds reveals that, if $\rho$ is a monotone increasing or decreasing function, the bounds do coincide, leading to
\begin{corollary}\label{maincor}
If in addition to the conditions of Theorem~\ref{maintheo} we assume that the ratio $\rho$ is monotone increasing or decreasing, then
$$
d_{\mathcal{W}}(P_1,P_2)= \frac{{\rm E}[\tau_1(\Theta_1;x)|\rho'(\Theta_1)|]}{{\rm E}[\rho(\Theta_1)]}.
$$
\end{corollary}

We now briefly discuss the relevance of the bounds from Theorem~\ref{maintheo}. The impact of a prior can best be measured by considering its effect on the resulting posterior distribution, and our result precisely allows us to quantify the difference in terms of Wasserstein distance between two posteriors resulting from distinct priors. As can be appreciated from~\eqref{bounds}, our bounds are sharp since they contain the same quantities in both the upper and lower bounds; this fact is further underlined by the equality in Corollary~\ref{maincor}.  The special case of one prior being the improper prior $p_1(\theta)=1$ allows us to retrieve the measure of the influence of the prior given in \citet{LRS17a}.

\section{Comparison of various priors for various distributions}\label{sec3}

We now illustrate the strength of Theorem~\ref{maintheo} by comparing popular choices of priors for parameters of three famous distributions, namely the normal, the Binomial and the Poisson models. We will work out in detail the case of the normal distribution. As we shall see, the bounds that we get allow us to conclude that, in all these cases, the difference between the resulting priors vanishes at the speed of $n^{-1}$ independently of the observations, but that at finite sample size $n$ the observations do play a role.

\subsection{Priors for the scale parameter of the normal distribution}\label{sec:nor}
Our first example concerns the normal $N(\mu,\sigma^2)$ distribution with probability density function
\begin{equation*}
x\mapsto\frac{1}{\sqrt{2\pi\sigma^2}} \exp \left(- \frac{(x-\mu)^2}{2\sigma^2} \right), x\in\R,
\end{equation*}
with location parameter $\mu\in\R$ and dispersion parameter $\sigma^2>0$. Our parameter of interest here is $\sigma^2$, and we consider $\mu$ to be fixed. Our motivation for studying various priors for $\sigma^2$ comes from \citet{KKF06} who consider a storm depth multiplier model to represent rainfall uncertainty, where the errors appear under multiplicative form and are assumed to be normal. They fix the mean $\mu$ but state that ``less is understood about the degree of rainfall uncertainty, i.e., the multiplier variance'' and therefore study various priors for $\sigma^2$. The two priors that they investigate are Jeffreys' prior and an Inverse Gamma prior, and we shall here show how our methods allow to directly measure the difference in impacts between these two priors.

Consider the  setting where $x = (x_1, \ldots , x_n)$ is a random sample from the $N(\mu, \sigma^2)$ population with fixed location $\mu$. The likelihood function $\ell(x;\sigma^2)$ of the normal model can be factorized into 
$$
\ell(x;\sigma^2) = (2\pi\sigma^2)^{-\frac{n}{2}} \exp \Bigg\{-\frac{1}{2} \sum_{i=1}^{n} {\frac{(x_i-\mu)^2}{\sigma^2}} \Bigg\}.
$$
The Jeffreys' prior is a popular prior as it is invariant under reparametrization and is proportional to the square-root of the Fisher information associated with the parameter of interest.  In the present setting it is proportional to $\frac{1}{\sigma^2}$ which leads to the posterior density being proportional to
$$
\sigma^2\mapsto(\sigma^2)^{-\frac{n}{2}-1}\exp \Bigg\{-\frac{1}{2} \sum_{i=1}^{n} {\frac{(x_i-\mu)^2}{\sigma^2}} \Bigg\}.
$$
The latter density can be recognized to be an Inverse Gamma distribution with  parameters  $(\frac{n}{2}, \frac{1}{2} \sum_{i=1}^{n}(x_i-\mu)^2)$. Now, it is easy to see that the Inverse Gamma prior 
$$
\sigma^2\mapsto \frac{\beta^\alpha}{\Gamma(\alpha)} (\sigma^2)^ {-\alpha-1} \exp \Bigg\{ -\frac{\beta}{\sigma^2} \Bigg\}
$$
with positive real parameters $(\alpha , \beta)$ yields, in combination with the likelihood, as posterior  also an Inverse Gamma distribution with parameters $(\frac{n}{2}+\alpha, \frac{1}{2}\sum_{i=1}^n(x_i-\mu)^2 + \beta)$. We consider Jeffreys' prior as first prior $P_1$ and the Inverse Gamma prior as second prior $P_2$, leading to the ratio
$$
\rho (\sigma^2) = \frac{p_2(\sigma^2)}{p_1(\sigma^2)} \propto \frac{(\sigma^2)^{-\alpha-1} \exp\left(-\frac{\beta}{\sigma^2}\right)}{(\sigma^2)^{-1}} =  (\sigma^2)^{-\alpha} \exp\left(-\frac{\beta}{\sigma^2}\right).
$$
With this in hand, Conditions (ii) and (iii)  of Theorem~\ref{maintheo} are readily fulfilled. Condition (i) is equivalent to requiring that $n/2+\alpha-1>0$ which is trivially the case as soon as we have at least two observations. Therefore, denoting by $\mu_1$ and $\mu_2$  the means of $P_1$ and $P_2$,  respectively, we can calculate the lower bound as follows:
\begin{eqnarray*}
d_{\mathcal W}(P_1, P_2) &\geq& |\mu_1-\mu_2|\\
 &=& \left|  \frac{\frac{1}{2} \sum_{i=1}^{n}(x_i-\mu)^2 }{ \frac{n}{2} - 1} -\frac{\frac{1}{2} \sum_{i=1}^n(x_i-\mu)^2 +\beta}{\frac{n}{2} + \alpha -1} \right| \\
&=& \frac{\left|  \frac{\alpha}{2} \sum_{i=1}^n(x_i-\mu)^2-(\frac{n}{2} - 1) \beta\right|}{(\frac{n}{2} + \alpha -1)(\frac{n}{2} - 1)}. 
\end{eqnarray*}
In order to obtain the upper bound, we need to calculate
\begin{eqnarray*}
\rho'(\sigma^2) &\propto&   (-\alpha) (\sigma^2)^{-\alpha-1}  \exp\left(-\frac{\beta}{\sigma^2}\right) + \frac{\beta}{(\sigma^2)^2} (\sigma^2)^{-\alpha} \exp\left(-\frac{\beta}{\sigma^2}\right)   \\
&=&  (\sigma^2)^{-\alpha-2} \exp\left(-\frac{\beta}{\sigma^2}\right) [-\alpha \sigma^2+ \beta ]
\end{eqnarray*}
and, writing $\Theta_1$ the random variable associated with $P_1$,
\begin{eqnarray*}
{\rm E} [\rho(\Theta_1)] &\propto& \int_{0}^{\infty}  (\sigma^2)^{-\alpha} \exp\left(-\frac{\beta}{\sigma^2}\right)\frac{ (\frac{1}{2} \sum_{i=1}^{n} (x_i-\mu)^2)^{\frac{n}{2}}}{\Gamma(\frac{n}{2})} (\sigma^2)^{-\frac{n}{2} -1}  \\
&\times& \exp{\left(-\frac{1}{\sigma^2} \frac{1}{2} \sum_{i=1}^{n} (x_i-\mu)^2\right)} d\sigma^2  \\
\small
&=&  \frac{ (\frac{1}{2} \sum_{i=1}^{n} (x_i-\mu)^2)^{\frac{n}{2}}}{\Gamma(\frac{n}{2})} \int_{0}^{\infty} (\sigma^2)^{-(\frac{n}{2} +\alpha) -1} \exp{\left(-\frac{1}{\sigma^2} \left(\beta + \frac{1}{2} \sum_{i=1}^{n} (x_i-\mu)^2\right)\right)} d\sigma^2 \\
\normalsize
&=& \frac{ (\frac{1}{2} \sum_{i=1}^{n} (x_i-\mu)^2)^{\frac{n}{2}}}{\Gamma(\frac{n}{2})}  \times \frac{\Gamma(\frac{n}{2}+\alpha)}{ (\beta + \frac{1}{2} \sum_{i=1}^{n} (x_i-\mu)^2)^{\frac{n}{2} + \alpha}}.
\end{eqnarray*}
The Stein kernel for the Inverse Gamma distribution with parameters $(\frac{n}{2}, \frac{1}{2} \sum_{i=1}^{n}(x_i-\mu)^2)$ corresponds to $\tau_1(\sigma^2) = \frac{(\sigma^2)^2}{\frac{n}{2}-1}$. Therefore we have
\begin{eqnarray*}
{\rm E} [\tau_1(\Theta_1) | \rho'(\Theta_1) |] &\propto&  \frac{ 1}{ (\frac{n}{2}-1)} {\rm E} \left[  (\Theta_1)^{-\alpha}  \exp\left(-\frac{\beta}{\Theta_1}\right) \left| -\alpha \Theta_1 + \beta \right|    \right]  \\
&\leq&  \frac{1}{(\frac{n}{2}-1)} \Bigg\{  \int_{0}^{\infty} \alpha (\sigma^2)^{-\alpha+1}  \exp\left(-\frac{\beta}{\sigma^2}\right) \frac{(\frac{1}{2} {\sum_{i=1}^{n} (x_i-\mu)^2)}^{\frac{n}{2}}}{\Gamma(\frac{n}{2})} \\
&\times& (\sigma^2)^{-\frac{n}{2} -1} \exp{\left(-\frac{1}{2\sigma^2} \sum_{i=1}^{n} (x_i-\mu)^2\right)} d\sigma^2  \\
&+& \int_{0}^{\infty} \beta (\sigma^2)^{-\alpha} \exp\left(-\frac{\beta}{\sigma^2}\right)  \frac{{(\frac{1}{2} \sum_{i=1}^{n} (x_i-\mu)^2})^{\frac{n}{2}}}{\Gamma(\frac{n}{2})} (\sigma^2)^{-\frac{n}{2}-1} \\
&\times&  \exp{\left(-\frac{1}{2\sigma^2} \sum_{i=1}^{n} (x_i-\mu)^2\right)}d\sigma^2  \Bigg\} \\
&=& \frac{ (\frac{1}{2} \sum_{i=1}^{n} (x_i -\mu)^2)^{\frac{n}{2}}}{(\frac{n}{2}-1)\Gamma(\frac{n}{2})} \Bigg\{ \int_{0}^{\infty} \alpha (\sigma^2)^{-(\alpha+\frac{n}{2})}  \\&\times& \exp{\left(-\frac{1}{\sigma^2} \left(\beta +\frac{1}{2} \sum_{i=1}^{n} (x_i - \mu)^2\right)\right) d\sigma^2} \\
&+& \beta\int_{0}^{\infty} (\sigma^2)^{-(\alpha+\frac{n}{2})-1} \exp{\left(-\frac{1}{\sigma^2} \left(\beta +\frac{1}{2} \sum_{i=1}^{n} (x_i - \mu)^2\right)\right)} d\sigma^2 \Bigg\} \\
&=& \frac{ (\frac{1}{2} \sum_{i=1}^{n} (x_i-\mu)^2)^{\frac{n}{2}}}{ (\frac{n}{2}-1) \Gamma(\frac{n}{2})} \Bigg[   \alpha \frac{\Gamma(\alpha+ \frac{n}{2} -1)}{(\beta+\frac{1}{2} \sum_{i=1}^{n} (x_i-\mu)^2)^{\alpha+ \frac{n}{2} -1}}   \\
&+& \beta \frac{\Gamma(\alpha+\frac{n}{2})}{(\beta+\frac{1}{2} \sum_{i=1}^{n} (x_i-\mu)^2)^{\alpha+ \frac{n}{2}}}    \Bigg] \\
&=& \frac{ (\frac{1}{2} \sum_{i=1}^{n} (x_i-\mu)^2)^{\frac{n}{2}} \Gamma(\alpha+\frac{n}{2} -1)}{ (\frac{n}{2}-1) \Gamma(\frac{n}{2}) (\beta +\frac{1}{2} \sum_{i=1}^{n} (x_i-\mu)^2)^{\alpha+\frac{n}{2}-1}} \left[  \alpha + \beta \frac{\alpha+\frac{n}{2}-1}{\beta +\frac{1}{2} \sum_{i=1}^{n} (x_i-\mu)^2}     \right]. 
\end{eqnarray*}
Putting the ends together, we thus get as upper bound
\begin{eqnarray*}
\frac{{\rm E} [\tau_1(\Theta_1) | \rho'(\Theta_1) |] }{{\rm E} [\rho(\Theta_1)]}&=&\frac{(\beta + \frac{1}{2} \sum_{i=1}^{n} (x_i-\mu)^2)}{(\frac{n}{2}+\alpha-1)(\frac{n}{2}-1)}\Bigg[  \alpha + \beta \frac{\alpha+\frac{n}{2}-1}{\beta +\frac{1}{2} \sum_{i=1}^{n} (x_i-\mu)^2}     \Bigg]\\
&=&\frac{\left(\frac{1}{2} \sum_{i=1}^{n} (x_i-\mu)^2\alpha+\frac{n}{2}\beta+(2\alpha-1)\beta\right)}{(\frac{n}{2}+\alpha-1)(\frac{n}{2}-1)}.
\end{eqnarray*}
We attract here the reader's attention to the fact that the same proportionality constant appears on the numerator and denominator, which is why we we left it out from our calculations for the sake of readability.

Finally we find that the Wasserstein distance between the posteriors based on the Jeffreys' and the Inverse Gamma priors for $\sigma^2$ is bounded as follows:
\small{
$$
\frac{\left| \alpha \sum_{i=1}^{n}(x_i-\mu)^2-(\frac{n}{2} - 1) \beta\right|}{(\frac{n}{2} + \alpha -1)(\frac{n}{2} - 1)}\leq d_{\mathcal W}(P_1,P_2)\leq \frac{\left(\frac{1}{2} \sum_{i=1}^{n} (x_i-\mu)^2\alpha+\frac{n}{2}\beta+(2\alpha-1)\beta\right)}{(\frac{n}{2}+\alpha-1)(\frac{n}{2}-1)}.
$$
}
We see that both the lower and upper bounds are of the order $O(n^{-1})$, which is a quantification of the well-known fact that the impact of the choice of the prior wanes asymptotically. Clearly, when $\alpha=\beta=0$, the distance equals zero which is natural as then the two priors coincide. It is interesting to remark that the  observations are centered around $\mu$ in the bounds, as this variation is precisely what the parameter $\sigma^2$ measures. Finally we remark that in the calculation of ${\rm E} [\tau_1(\Theta_1) | \rho'(\Theta_1) |]$ we have used the triangular inequality for $\left| -\alpha \Theta_1 + \beta \right|$, while alternative upper boundings or direct calculations with absolute values would lead to a sharper upper bound. However, those will not be as readable as our bounds, which anyway are of the desired order  of~$n^{-1}$.

\subsection{{Priors for the success parameter in a Binomial model}}

Many random phenomena worth studying have binary outcomes and therefore can be modelled using the famous Binomial distribution $Bin (n,\theta)$ with probability mass function
$$
x\mapsto {{n} \choose {x}} \theta^x (1-\theta)^{n-x}
$$
where $x \in \{1, \ldots, n\}$ is the number of observed successes, the natural number $n$ denotes the number of binary experiments and $\theta \in (0,1)$ stands for the success parameter. In this part, we suppose $n$ to be fixed and we focus on the success parameter. \citet{Ghad18} studies the influence of various priors, including the Beta and Haldane priors, on the posterior by bounding the Wasserstein distance between each resulting posterior and the data-only posterior obtained via the uniform prior. We will show how to compare the Haldane and Beta priors directly by means of our Theorem~\ref{maintheo}.

The Haldane prior is an improper prior representing complete uncertainty  and its expression is given by
$$
p_1 (\theta) = \frac{1}{\theta (1-\theta)}.
$$
The resulting posterior density is of the form
$$
\theta\mapsto {{n} \choose {x}} \theta^{x-1} (1-\theta)^{n-x-1}
$$
which represents a Beta distribution with parameters $(x,n-x)$. In order to be a proper distribution it requires at least one success and one failure in the sequence of experiments. Choosing as second prior the conjugate Beta prior with parameters $\alpha,\beta>0$ yields the  $Beta(\alpha+x, n-x+\beta)$ posterior.  The ratio of prior densities corresponds here to 
$$
\rho(\theta) =\frac{p_2(\theta)}{p_1(\theta)}\propto \theta^{\alpha} (1-\theta)^{\beta}.
$$
One readily sees that the conditions of Theorem~\ref{maintheo} are satisfied. The lower bound thus is given by
$$
d_{\mathcal{W}}(P_1, P_2)  \geq \left| \frac{x}{n}-\frac{\alpha + x}{n+\alpha+\beta} \right|  =  \frac{\left|n\alpha - (\alpha+\beta)x\right|}{n(n+\alpha+\beta)}.
$$
The Stein kernel for the Beta distribution with parameters $(x,n-x)$ is $\tau_1(\theta) = \frac{\theta (1-\theta)}{n}$. Proceeding as in Section~\ref{sec:nor}, we get the upper bound
$$
d_{\mathcal{W}}(P_1, P_2) \leq \frac{1}{n} \left( \alpha + (\beta-\alpha) \frac{\alpha + x}{\alpha + \beta + n} \right)
$$
It can be seen both the lower and upper bound are of the order of $O(n^{-1})$. This rate of convergence remains even in the extreme cases $x=0$ and $x=n$. 

\subsection{{Priors for the Poisson model}}

As last example we consider the most popular count distribution, namely the Poisson distribution with probability mass function 
$$
x\mapsto\frac{\exp{(-\theta)} \theta^x}{x!}
$$
where $x\in N$ is the number of occurrences and the parameter $\lambda > 0$  indicates the average number of events in a given time interval. Classical choices of priors are an exponential prior, Jeffreys' prior or the conjugate Gamma prior. The latter was used in e.g. the study of asthma mortality rates by \citet{Gelman04}. We opt here to compare in all generality two Gamma priors with non-negative real parameters $(\alpha_1,\beta_1)$ and $(\alpha_2,\beta_2)$, respectively, as this contains as special cases the exponential prior ($\alpha=1$), the uniform prior ($\alpha=1,\beta=0$) and even the Jeffreys' prior which is proportional to $\theta^{-1/2}$ ($\alpha=1/2,\beta=0$).

The likelihood function $\ell(x;\theta)$ for a set of sampling data $x = (x_1,  \ldots, x_n)$ coming from a Poisson model reads 
$$
\ell(x;\theta) = \frac{\exp{(-n\theta)}\theta^{\sum_{i=1}^{n} x_i}}{\prod_{i=1}^{n} x_i!}.
$$
In combination with a conjugate Gamma prior $\frac{\beta_j^{\alpha_j}}{\Gamma(\alpha_j)}\theta^{\alpha_j-1}\exp(-\beta_j\theta)$ this yields a Gamma posterior distribution with updated parameters $(\sum_{i=1}^nx_i+\alpha_j,n+\beta_j)$, $j=1,2$. Now, a quantity of special interest is the ratio 
$$
\rho(\theta)=\frac{p_2(\theta)}{p_1(\theta)}\propto \theta^{\alpha_2-\alpha_1}\exp\left(-(\beta_2-\beta_1)\theta\right).
$$
While it nicely satisfies the conditions of Theorem~\ref{maintheo}, certain combinations of the prior parameters lead to a monotone increasing or decreasing behavior, and hence the applicability of Corollary~\ref{maincor}. Simple manipulations reveal that these parameter combinations are given by $\alpha_1<\alpha_2\cap\beta_1>\beta_2$ (increasing) and $\alpha_1>\alpha_2\cap\beta_1<\beta_2$ (decreasing), where in each case one of the strict inequalities can  turn into a potential equality. The Stein kernel of the Gamma distribution with parameters $(\sum_{i=1}^nx_i+\alpha_1,n+\beta_1)$ being $\tau_1(\theta)=\frac{\theta}{n+\beta_1}$, straightforward calculations lead to
$$
{\rm E}[\rho(\Theta_1)]\propto\frac{\Gamma(\sum_{i=1}^nx_i+\alpha_2)}{(n+\beta_2)^{\sum_{i=1}^nx_i+\alpha_2}}
$$
and 
$$
{\rm E}[\rho'(\Theta_1)\tau_1(\Theta_1)]\propto\frac{(\alpha_2-\alpha_1)\frac{\Gamma(\sum_{i=1}^nx_i+\alpha_2)}{(n+\beta_2)^{\sum_{i=1}^nx_i+\alpha_2}}-(\beta_2-\beta_1)\frac{\Gamma(\sum_{i=1}^nx_i+\alpha_2+1)}{(n+\beta_2)^{\sum_{i=1}^nx_i+\alpha_2+1}}}{n+\beta_2},
$$
and hence
$$
d_{\mathcal{W}}(P_1, P_2)=\frac{1}{n+\beta_2}\left(\alpha_2-\alpha_1-(\beta_2-\beta_1)\frac{\sum_{i=1}^nx_i+\alpha_2}{n+\beta_2}\right).
$$
This exact distance is of the desired order of $O(n^{-1})$. In case of parameter combinations that lead to no monotone behavior of the ratio $\rho$ we can  build upper and lower bounds as done in the previous sections; we leave this as exercise for the reader.

\

\noindent ACKNOWLEDGMENTS:

\noindent This research is supported by a BOF Starting Grant of Ghent University.

%
%
%
%

\bibliographystyle{elsarticle-num-names}

\bibliography{GL18SPL}

\end{document}